\documentclass{article}

\usepackage{arxiv}

\usepackage[utf8]{inputenc} 
\usepackage[T1]{fontenc}    
\usepackage{hyperref}       
\usepackage{url}            
\usepackage{booktabs}       
\usepackage{amsfonts}       
\usepackage{nicefrac}       
\usepackage{microtype}      
\usepackage{lipsum}

\usepackage{graphicx}

\title{Exponential Fourth Order Schemes for Direct Zakharov-Shabat problem}

\author{
  Sergey Medvedev$^{1,2,*}$, Irina Vaseva$^{1,2}$, Igor
Chekhovskoy$^{2,1}$, Mikhail Fedoruk$^{2,1}$\\
$^{1}$ Institute of Computational Technologies, SB RAS, Novosibirsk
630090, Russia,\\
$^{2}$ Novosibirsk State University, Novosibirsk 630090, Russia,\\
* Corresponding author: medvedev@ict.nsc.ru
}

\begin{document}
\maketitle

\begin{abstract}
We propose two finite-difference algorithms of fourth order of accuracy for solving the initial problem of the Zakharov-Shabat system. Both schemes have the exponential form and conserve quadratic invariant of Zakharov-Shabat system. The second scheme contains the spectral parameter in exponent only and allows to apply the fast computational algorithm.
\end{abstract}


\section{Introduction}
Recently, there has been great interest in the so-called
nonlinear Fourier transform (NFT), which is a generalization
of the ordinary linear Fourier transform. The NFT makes it
possible to find exact solutions for nonlinear equations, such
as the nonlinear Schr\"odinger equation (NLSE) and the
Korteweg de Vries (KdV) equation, interpreting the evolution
of optical fields by analogy with the linear Fourier transform
as an evolution of a certain set of frequencies.
For the first time, this idea for the NLSE was proposed by Zakharov and Shabat in
1971~\cite{ZakharovShabat1972}. They showed that the
NLSE can be integrated by the inverse scattering problem
(IST) method, previously applied to the Korteweg de Vries
(KdV) equation. The NLSE describes the envelope for wave
beams, therefore it is used in many areas of physics where
there are wave systems.
Despite a large number of articles~\cite{Yousefi2014II, Turitsyn2017Optica, Vasylchenkova2019a} devoted to NFT, the
development of the accurate and fast numerical algorithms for
NFT still remains an actual mathematical problem.

In this paper, we consider the direct spectral problem which consists in the numerical solution of the Zakharov-Shabat (ZS) system. The integration of this system is the first step in the general scheme of NFT. In addition, we limit ourselves to constructing one-step finite-difference schemes for solving the ZS system and do not touch the question of methods for finding discrete spectral parameters and phase coefficients for them.

We present the general necessary conditions for the transition operator for fourth-order one-step difference schemes for linear systems of first-order differential equations. Then we give two examples of such schemes. Both schemes are the exponential fourth order ones, and we show their connection with the Magnus decomposition. The main property of such schemes is to conserve the quadratic invariant for continuous spectral parameters. The second scheme contains the spectral parameter only in the exponent. This property allows one to apply a fast algorithm for the numerical solution of the ZS system.

The final part of the article presents comparisons of numerical computations using the proposed schemes and other well-known schemes: the Boffetta-Osborne second-order scheme~\cite{Boffetta1992a}, Runge-Kutta fourth-order scheme and the fourth-order conservative scheme~\cite{Medvedev2019_OL}, which is a generalization of the Boffetta-Osborn scheme.

\section{The direct Zakharov-Shabat problem}
The standard NLSE is a basic model for the pulse propagation along an ideally lossless and noiseless fiber
\begin{equation}\label{nlse}
i\frac{\partial q}{\partial z}+\frac{\sigma}{2}\frac{\partial^2 q}{\partial t^2}+|q|^2q=0,
\end{equation}
where $q=q(t,z)$ is a slow-varying complex optical field envelope, the variable $z$ is the distance along the optical fiber, $t$ is a time variable; $\sigma =-1$ and $\sigma =1$ corresponds to the normal and anomalous dispersion, respectively. Eq.~(\ref{nlse}) is written in the moving coordinate system and describes the propagation of pulses $q(t,z)$ in optical fibers. Therefore, the initial data are almost stationary, and the Cauchy problem is solved with the initial conditions as follows:
$$\left.q(t,z)\right|_{z=z_0}=q_0(t).$$

The mathematical method suggested by Zakharov and Shabat~\cite{ZakharovShabat1972} allows to integrate the NLSE. The method, widely known as the Nonlinear Fourier Transform (NFT), allows transforming signal into nonlinear Fourier spectrum, which defined by the solution of the Zakharov-Shabat problem (ZSP).

The equation (\ref{nlse}) can be written as a condition of compatibility
\begin{equation}
\frac{\partial L}{\partial z}=ML-LM
\end{equation}
of two linear equations
\begin{equation}\label{LM}
L\Psi=\zeta\Psi,\quad \frac{\partial\Psi}{\partial z}=M\Psi,
\end{equation}
where $\Psi(t)$ is a complex vector function of a real argument $t$,
\begin{equation}
L=i\left(\begin{array}{cc}\partial_t&-q\\-\sigma q^*&-\partial_t\end{array}\right),\quad
M=i\left(\begin{array}{cc}
\sigma\frac{\partial^2}{\partial t^2}+\frac{1}{2}|q|^2&
-\sigma q\frac{\partial}{\partial t}-\frac{1}{2}\sigma\frac{\partial q}{\partial t}\\
-q^*\frac{\partial}{\partial t}-\frac{1}{2}\frac{\partial q^*}{\partial t}&
-\sigma\frac{\partial^2}{\partial t^2}-\frac{1}{2}|q|^2\end{array}\right).
\end{equation}

The first equation in (\ref{LM}) is the eigenvalue problem for the operator $L$.
For $\sigma=-1$ the operator $L$ is Hermitian ($L=L^\dagger\equiv (L^*)^T$), therefore the complex spectral parameter $\zeta=\xi+i\eta$ becomes real $\zeta=\xi\in \mathbb{R}$. There is no such restriction for $\sigma=1$. In this case the problem has the continuous and discrete spectra. The continuous spectrum lies on the real axis and the discrete spectrum is in the upper half plane $\mbox{Im}(\zeta) > 0$ .

Also the first equation in (\ref{LM}) can be rewritten as an evolutionary system
\begin{equation}\label{psit}
\frac{d \Psi(t)}{dt}=Q(t){\Psi}(t),
\end{equation}
where $q=q(t,z)$ and
$${\Psi}(t)=\left(\begin{array}{c}\psi_1(t)\\\psi_2(t)\end{array}\right),\quad
Q(t)=\left(\begin{array}{cc}-i\zeta&q\\-\sigma q^*&i\zeta\end{array}\right).
$$
Here $z$ is a parameter, that we will skip further.

For $\sigma=1$ the matrix $Q$ is skew-Hermitian $Q=-Q^\dagger$. This leads to conservation of the quadratic energy invariant for the real spectral parameters. Moreover, the system~(\ref{psit}) can be written in a gradient form as follows:
\begin{equation}
\left(\begin{array}{c}\psi_1\\\psi_2\end{array}\right)_t=
J\left(\begin{array}{c}\psi_1\\\sigma\psi_2\end{array}\right)=
J\left(\begin{array}{c}\frac{\partial H}{\partial\psi_1^*}\\\frac{\partial H}{\partial\psi_2^*}\end{array}\right),
\quad J=\left(\begin{array}{cc}-i\zeta&\sigma q\\-\sigma q^*&i\sigma\zeta\end{array}\right)
\end{equation}
where $H=|\psi_1|^2+\sigma|\psi_2|^2$. For real spectral parameters $\zeta=\xi$ the matrix $J$
is skew-Hermitian for any $\sigma=\pm 1$ and, consequently, $H$ conserves. The value $H$ and the matrix $Q$ can be written using Pauli matrices $\sigma_0$ and $\sigma_3$ as follows:
\begin{equation}\label{HQ}
H=\left\{\begin{array}{ccc}
({\Psi}^*,\sigma_0{\Psi}),&\mbox{for }&\sigma=1\\
({\Psi}^*,\sigma_3{\Psi}),&\mbox{for }&\sigma=-1
\end{array}\right.,\quad Q=\left\{\begin{array}{ccc}
J\sigma_0,&\mbox{for }&\sigma=1\\
J\sigma_3,&\mbox{for }&\sigma=-1
\end{array}\right..
\end{equation}

Assuming that $q(t)$ decays rapidly when $t\to \pm \infty $, the specific solutions (Jost functions) for ZSP (\ref{psit}) can be derived as:
\begin{equation}\label{psi0}
\Psi =\left(
\begin{array}{c}
\psi_{1}\\\psi_{2}
\end{array}
\right) = \left(
\begin{array}{c}
e^{-i\zeta t}\\0
\end{array}
\right)[1+o(1)],\quad t\to-\infty,
\end{equation}
and
\begin{equation}\label{psi0right}
\Phi =\left(
\begin{array}{c}
\phi_{1}\\\phi_{2}
\end{array}
\right) = \left(
\begin{array}{c}
0\\e^{i\zeta t}
\end{array}
\right)[1+o(1)],\quad t\to \infty,
\end{equation}
Then we obtain the Jost scattering coefficients $a(\xi )$ and $b(\xi )$ as follows:
\begin{equation}\label{ab}
a(\xi)=\lim_{t\to\infty}\,\psi_1(t,\xi)\,e^{i\xi t},\quad b(\xi)=\lim_{t\to\infty}\,\psi_2(t,\xi)\,e^{-i\xi t}.
\end{equation}
The functions $a(\xi )$ and $b(\xi )$ can be extended to the upper half-plane $\xi \to \zeta $, where $\zeta $ is a complex number with the positive imaginary part $\eta = \mbox{Im}\zeta >0$. The spectral data of ZSP (\ref{psit}) are determined by $a(\zeta )$ and $b(\zeta )$ in the following way:\\
(1) $K$ zeros of $a(\zeta)=0$ define the discrete spectrum $\{\zeta_k\}$, $k=\overline{0,K-1}$ of ZSP~(\ref{psit}) and phase coefficients
$$r_k=\left.\frac{b(\zeta)}{a'(\zeta)}\right|_{\zeta=\zeta_k},\quad\mbox{where}\quad a'(\zeta)=\frac{da(\zeta)}{d\zeta};$$
(2) the continuous spectrum is determined by the reflection coefficient
$r(\xi)=b(\xi)/a(\xi)$, $\xi\in\mathbb{R}$.

These spectral data were defined using the "left" boundary condition~(\ref{psi0}). Both conditions~(\ref{psi0}) and~(\ref{psi0right}) can be used to calculate the coefficient $b(\zeta_{k} )$ of the discrete spectrum:
\begin{equation}\label{bidir}
\Psi(t,\zeta_k)=\Phi(t,\zeta_k)b(\zeta_k).
\end{equation}

For real values of the spectral parameter $\zeta=\xi$ we have invariant $H$.
Taking into account the boundary conditions~(\ref{psi0}), we get $H=1$.

In addition, the following trace formula is valid~\cite{Ablowitz1981}:
\begin{equation}
C_n=-\frac{1}{\pi}\int\limits_{-\infty}^\infty\,(2i\xi)^n\,\ln|a(\xi)|^2\,d\xi+
\sum\limits_{k=0}^{K-1}\,\frac{1}{n+1}
\left[(2i\zeta_k^*)^{n+1}-(2i\zeta_k)^{n+1}\right],
\end{equation}
which connects the NLSE integrals $C_n$ with the coefficient $a(\xi)$ and the discrete spectrum $\zeta_k$. The first integrals have the form
$$C_0=\int\limits_{-\infty}^\infty|q|^2dt,\enskip C_1=\int\limits_{-\infty}^\infty qq^*_tdt,\enskip C_2=\int\limits_{-\infty}^\infty(qq^*_{tt}+|q|^4)dt,\enskip
C_3=\int\limits_{-\infty}^\infty(qq^*_{ttt}+4|q|^2qq^*_t+|q|^2q^*q_t)dt.$$
This formula with $n=0$
\begin{equation} \label{Parseval}
C_0=-\frac{1}{\pi } \int _{-\infty }^{\infty } \ln |a(\xi )|^{2} d\xi +\sum _{k=0}^{K-1} \left[2i\left(\zeta _{k}^{*} -\zeta _{k} \right)\right]
\end{equation}
is called the Parseval nonlinear equality and is used to verify the numerical calculations and the consistency of the continuous and discrete spectra found. The first term on the right-hand side of Eq.~(\ref{Parseval}) refers to the continuous spectrum energy:
\begin{equation} \label{Econt}
E_{c} =-\frac{1}{\pi } \int _{-\infty }^{\infty } \ln |a(\xi )|^{2} d\xi .
\end{equation}

We solve the system~(\ref{psit}). The matrix~$Q(t)$ linearly depends on the complex function~$q(t)$ which is given in the whole nodes of the uniform grid with a step~$\tau$ on the interval $[-L,L]$. Let us note main features of the discrete problem:
\begin{itemize}
	\item Since the matrix $Q$ is defined on a uniform grid, the unknown function $\Psi $ must also be computed on a uniform grid with the same step. Therefore, the Runge-Kutta methods (RK) cannot be used on such grid. If, for example, we use an explicit 4th order RK scheme, then we need to take the computational grid with a double step $2\tau$. In this case, values of $Q_n$ will be used unequally.
	
	\item For small values of the potential $|q(t)|<<|\zeta|$ and $\mbox{Im}\,\zeta>0$, ZSP has exponentially growing and decreasing solutions, thus A-stability of finite-difference methods is required~\cite{Dahlquist1963}. The method is called A-stable if all solutions of the equation $\partial x/\partial t=\lambda x$ tend to zero at $\mbox{Re}\,\lambda<0$ and fixed step $h$. The second barrier of Dahlquist restricts the use of multi-step methods~\cite{Hairer1987}. It means that there are no explicit A-stable multishep methods for the Eq.~(\ref{psit}), and the 2nd order of convergence is maximal for implicit multi-step methods.
	
	\item The ZSP has a second order matrix, therefore, the inverse matrices and the matrix exponential can be easily calculated. This allows us to include practically any functions of the matrix~$Q$ in the difference schemes.
	
	\item To calculate the spectral data, it is necessary to solve the ZSP for a large number of values~$\zeta$ at a fixed potential~$q(t)$. This should be taken into account when implementing algorithms.
	
\end{itemize}

\section{General theory of one-step schemes}

Let us consider the problem (\ref{psit}) in a general case. We need to solve the equation
\begin{equation}\label{xt}
Dx=Q(t)x,\quad D=\frac{d}{dt},
\end{equation}
where $x=x(t)\in \mathbb{C}^n$, using a one-step algorithm
\begin{equation}\label{Tx}
x_{n+1}=Tx_n,
\end{equation}
where $T$ is the transition operator, $x_n=x(t_n)$, $t_n=n\tau$, $\tau$ is a step of the uniform grid.

We differentiate Eq.~(\ref{xt}) and get the expressions for the derivatives $D^k x$ up to $5$-th order as follows:
\begin{equation}\label{d1}
\begin{array}{l}
Dx=Qx,\\[2mm]
D^2x=(DQ)x+QDx,\\[2mm]
D^3x=(D^2Q)x+2(DQ)(Dx)+QD^2x,\\[2mm]
D^4x=(D^3Q)x+3(D^2Q)(Dx)+3(DQ)(D^2x)+QD^3x,\\[2mm]
D^5x=(D^4Q)x+4(D^3Q)(Dx)+6(D^2Q)(D^2x)+4(DQ)(D^3x)+QD^4x.
\end{array}
\end{equation}
Let us introduce the notation for the right-hand side of Eq.~(\ref{d1}) and derivatives $Q^{(k)} = D^kQ$
\begin{equation}\label{d2}
D^kx=Q_{k}x,\quad Q_1=Q.
\end{equation}
Using (\ref{d1}) and (\ref{d2}) we find recurrence relations for $Q_k$ as follows:
$$
\begin{array}{l}
Q_2=Q^{(1)}+Q^2,\\[2mm]
Q_3=Q^{(2)}+2Q^{(1)}Q+QQ^{(1)}+Q^3,\\[2mm]
Q_4=Q^{(3)}+3Q^{(2)}Q+QQ^{(2)}+3(Q^{(1)})^2+3Q^{(1)}Q^2+2QQ^{(1)}Q+Q^2Q^{(1)}+Q^4,\\[2mm]
Q_5=Q^{(4)}+4Q^{(3)}Q+QQ^{(3)}+6Q^{(2)}Q^{(1)}+4Q^{(1)}Q^{(2)}+6Q^{(2)}Q^2+3QQ^{(2)}Q+Q^2Q^{(2)}+\\[2mm]
+8(Q^{(1)})^2Q+4Q^{(1)}QQ^{(1)}+3Q(Q^{(1)})^2+4Q^{(1)}Q^3+3QQ^{(1)}Q^2+2Q^2Q^{(1)}Q+Q^3Q^{(1)}+Q^5.
\end{array}
$$

Let us derive the Taylor series of $x(t)$ at the point $t$, such that $t_n=t+\overline{s}\tau$, $t_{n+1}=t+s\tau$, $\overline{s} =s- 1$:
\begin{equation}\label{left}
x(t_{n+1})= x+s\tau Dx+\frac{(s\tau)^2}{2!}D^2x+\frac{(s\tau)^3}{3!}D^3x+\frac{(s\tau)^4}{4!}D^4x+\frac{(s\tau)^5}{5!}D^5x+O(\tau^6),
\end{equation}
\begin{equation}\label{right}
x(t_n) =
\displaystyle
x+\overline{s}\tau Dx+\frac{(\overline{s}\tau)^2}{2!}D^2x+\frac{(\overline{s}\tau)^3}{3!}D^3x+
+\frac{(\overline{s}\tau)^4}{4!}D^4x +\frac{(\overline{s}\tau)^5}{5!}D^5x+O(\tau^6).
\end{equation}
Then we denote the terms of (\ref{left}) and (\ref{right}):
\begin{equation}\label{d3}
L_k=\frac{s^k}{k!}Q_{k},\quad R_k=\frac{\overline{s}^k}{k!}Q_{k}
\end{equation}
and write the expansion of (\ref{Tx}) up to $5$-th order
\begin{equation}
(E+L_1+L_2+L_3+L_4+L_5)=(T_0+T_1+T_2+T_3+T_4+T_5)(E+R_1+R_2+R_3+R_4+R_5).
\end{equation}
After equating the terms of the same order we get
\begin{equation}\label{t1}
L_1=R_1+T_1,
\end{equation}
\begin{equation}\label{t2}
L_2=R_2+T_1R_1+T_2,
\end{equation}
\begin{equation}\label{t3}
L_3=R_3+T_1R_2+T_2R_1+T_3,
\end{equation}
\begin{equation}\label{t4}
L_4=R_4+T_1R_3+T_2R_2+T_3R_1+T_4,
\end{equation}
\begin{equation}\label{t5}
L_5=R_5+T_1R_4+T_2R_3+T_3R_2+T_4R_1+T_5.
\end{equation}

Now we can derive recurrence relations for $T_k$. From Eq.~(\ref{t1}) we get
\begin{equation}
T_1=L_1-R_1=s Q-\overline{s}Q=Q.
\end{equation}
Therefore, for first-order approximation, the expansion of transition operator~$T$ must begin with $T\approx E+\tau Q$ and we need to know the value of $Q$ at the point $t$.

From Eq.~(\ref{t2}) we get
\begin{equation}
T_2=L_2-R_2-T_1R_1=\frac{s^2-\overline{s}^2}{2!}Q_2-\overline{s}Q^2=\frac{2s-1}{2!}Q_2-\overline{s}Q^2.
\end{equation}
To find $T_2$ we need to know the values of $Q^2$ and $Q^{(1)}$ at the point $t$. If we do not have the analytical expression of $Q(t)$, we need to know the value of $Q$ at two different points to use finite differences to calculate $Q^{(1)}$. Otherwise, we can set $s=1/2$ and zero the coefficient at $Q_2$. Then we only need to know the value of $Q^2$ at the point $t$.

From Eq. (\ref{t3}) we get
\begin{equation}\label{T3a}
\begin{array}{lll}
T_3 & = & \displaystyle L_3-R_3-T_1R_2-T_2R_1=\frac{s^3-\overline{s}^3}{3!}Q_3-Q\frac{\overline{s}^2}{2!}Q_2-\overline{s}T_2Q=\\[4mm]
& = & \displaystyle \frac{3s^2-3s+1}{3!}Q_3-\frac{\overline{s}^2}{2!}QQ_2-\frac{(2s-1)\overline{s}}{2!}Q_2Q+\overline{s}^2Q^3.
\end{array}
\end{equation}
Equation
$$s^3-\overline{s}^3=3s^2-3s+1=0$$
has no real roots, therefore we can not get rid of the term with $Q_3$ varying $s$. It means that to use any scheme of order higher then $2$ we must know $Q^{(2)}$ or values of $Q$ at three different points.

From Eq. (\ref{t4}) we get
\begin{equation}\label{T4a}
T_4 = \displaystyle
L_4-R_4-T_1R_3-T_2R_2-T_3R_1=
\frac{s^4-\overline{s}^4}{4!}Q_4-\frac{\overline{s}^3}{3!}QQ_3-\frac{\overline{s}^2}{2!}T_2Q_2-\overline{s}T_3Q.
\end{equation}
Since equation
$$s^4-\overline{s}^4=(2s-1)(2s^2-2s+1)=0$$
has only one real root $s=1/2$, this is the only way to get rid of $Q_4$, that contains $Q^{(3)}$.

From Eq.~(\ref{t5}) we get
\begin{equation}
T_5=L_5-R_5-T_1R_4-T_2R_3-T_3R_2-T_4R_1=\frac{s^5-\overline{s}^5}{5!}Q_4+...\,.
\end{equation}
Since equation
$$s^5-\overline{s}^5=0$$
has no real roots, we can not zero this coefficient varying $s$.

Thus we formulate

\textbf{Theorem.} Any one-step finite-difference scheme~(\ref{Tx}) approximates the equation (\ref{xt}) with a fourth order of accuracy if and only if the expansion of the transition operator~$T$ for the fixed~$s$ has a form
\begin{equation}\label{Tx4}
T=E+\tau Q+\tau^2 T_2+\tau^3T_3+\tau^4T_4+O(\tau^5),
\end{equation}
where
\begin{equation}\label{T2a}
T_2=\frac{2s-1}{2!}Q_2-\overline{s}Q^2,
\end{equation}
\begin{equation}\label{T3b}
T_3=\frac{3s^2-3s+1}{3!}Q_3-\frac{\overline{s}^2}{2!}QQ_2-\frac{(2s-1)\overline{s}}{2!}Q_2Q+\overline{s}^2Q^3,
\end{equation}
\begin{equation}\label{T4b}
T_4=\frac{(2s-1)(2s^2-2s+1)}{4!}Q_4-\frac{\overline{s}^3}{3!}QQ_3-\frac{\overline{s}^2}{2!}T_2Q_2-\overline{s}T_3Q
\end{equation}
and the coefficients~$Q_k$ are expressed through the matrix $Q$ and its derivatives
\begin{equation}
Q_2=Q^{(1)}+Q^2,
\end{equation}
\begin{equation}
Q_3=Q^{(2)}+2Q^{(1)}Q+QQ^{(1)}+Q^3,
\end{equation}
\begin{equation}
Q_4=Q^{(3)}+3Q^{(2)}Q+QQ^{(2)}+3(Q^{(1)})^2+3Q^{(1)}Q^2+2QQ^{(1)}Q+Q^2Q^{(1)}+Q^4.
\end{equation}

\section{Examples of schemes}
Let us consider examples of constructing fourth-order schemes that satisfy the conditions of the theorem.

\subsection{Constant matrix $Q$}
\textbf{Corollary 1.}
If the matrix~$Q$ is constant, then $Q_n=Q^n$ and the expansion of the matrix~$T$ does not depend on~$s$ and has the form
\begin{equation}\label{Texp}
T=E+\tau Q+\frac{\tau^2}{2!}Q^2+\frac{\tau^3}{3!}Q^3+\frac{\tau^4}{4!}Q^4+O(\tau^5).
\end{equation}
It is clear that Eq. (\ref{Texp}) is the expansion of the matrix exponential~$\exp(\tau Q)$. Since it is an exact solution for the system with a constant matrix, then the one-step scheme with the exponential form has the infinity order of approximation for transition operator~$T=\exp(\tau Q)$.

\subsection{Symmetrical case}
As mentioned before, the second order scheme does not contain the derivative $Q^{(1)}$ if and only if $s=1/2$. It means that the transition matrix $T$ depends only on $Q(t_n+\tau/2)$. Such choice of $s$ corresponds to the center of the interval and will be called the symmetrical case. The expansion of the matrix~$T$ for the fourth order schemes in the symmetrical case gets rid of some terms and has only the dependence on $Q^{(1)}$ and $Q^{(2)}$. Therefore, it is necessary to use $Q$ at least at three points for the fourth order scheme.

\textbf{Corollary 2.}
The expansion (\ref{Tx4}) of the matrix $T$ in the symmetrical case has the form as follows:
\begin{equation}\label{T4se1}
\begin{array}{l}
T = \displaystyle E+\tau Q+\frac{1}{2}\tau^2Q^2+\frac{\tau^3}{3!}Q^3+\frac{\tau^3}{24}Q^{(2)}+\frac{\tau^3}{12}\left(Q^{(1)}Q-QQ^{(1)}\right)+\\
\quad \displaystyle +\frac{\tau^4}{4!}Q^4+\frac{\tau^4}{48}\left(QQ^{(2)}+Q^{(2)}Q\right)+\frac{\tau^4}{24}\left(Q^{(1)}Q^2-Q^2Q^{(1)}\right).
\end{array}
\end{equation}
Approximating the derivatives in Eq.~(\ref{T4se1}) by central finite differences of the second order
\begin{equation}\label{Q12}
Q^{(1)}_{n+\frac{1}{2}}=\frac{Q_{n+\frac{3}{2}}-Q_{n-\frac{1}{2}}}{2\tau}+O(\tau^2),\quad Q^{(2)}_{n+\frac{1}{2}}=\frac{Q_{n+\frac{3}{2}}-2Q_{n+\frac{1}{2}}+Q_{n-\frac{1}{2}}}{\tau^2}+O(\tau^2).
\end{equation}
we retain the fourth order of accuracy of the operator
\begin{equation}\label{T4se1r}
\begin{array}{ll}
T_{n+\frac{1}{2}} & = \displaystyle E+\tau Q_{n+\frac{1}{2}}+\frac{\tau^2}{2}Q_{n+\frac{1}{2}}^2+\frac{\tau^3}{3!}Q^3+
\frac{\tau^4}{4!}Q_{n+\frac{1}{2}}^4+\\
& \displaystyle
+\frac{\tau^3}{12}\left(Q^{(1)}_{n+\frac{1}{2}}Q_{n+\frac{1}{2}}-Q_{n+\frac{1}{2}}Q^{(1)}_{n+\frac{1}{2}}\right)+
\frac{\tau^3}{24}Q^{(2)}_{n+\frac{1}{2}}+\\
& \displaystyle
+\frac{\tau^4}{48}\left(Q_{n+\frac{1}{2}}Q^{(2)}_{n+\frac{1}{2}}+Q^{(2)}_{n+\frac{1}{2}}Q_{n+\frac{1}{2}}\right)+
\frac{\tau^4}{24}\left(Q^{(1)}_{n+\frac{1}{2}}Q_{n+\frac{1}{2}}^2-Q_{n+\frac{1}{2}}^2Q^{(1)}_{n+\frac{1}{2}}\right)+O(\tau^5).
\end{array}
\end{equation}

\subsection{Exponential form}
The formula~(\ref{T4se1}) is an expansion of the exponent
\begin{equation}\label{exp4}
T=\exp\left\{\tau F_1+\tau^3 F_3\right\}+O(\tau^5),
\end{equation}
where
\begin{equation}\label{power4}
F_1=Q,\quad F_3=\frac{1}{24}Q^{(2)}+\frac{1}{12}\left(Q^{(1)}Q-QQ^{(1)}\right),\quad F_2=F_4=0.
\end{equation}
Replacing the derivatives by the finite differences (\ref{Q12}) we get the exponential scheme with the fourth order of accuracy (ES4). The scheme conserves energy for the real spectral parameters~$\zeta=\xi$.

Exponential form (\ref{exp4}) of the scheme allows approximating it by the Pad\'e approximation and apply to multidimensional systems.

\subsection{Magnus expansion}

Exponential form of Eq.~(\ref{exp4}) follows from the Magnus expansion~\cite{Magnus1954}. It provides an exponential representation of the exact evolution operator of the system (\ref{psit})
$$\Psi(t)=U(t,0)\Psi(0),\quad U(t,0)=e^{\Omega(t)},$$
which is constructed as a series expansion referred to as Magnus expansion:
$\Omega(t)=\sum\limits_{k=0}^\infty \Omega_k(t).$
First terms of this series have forms as follows:
$$\Omega_1(t)=\int\limits_{0}^{t}dt_1 Q(t_1),\quad
\Omega_2(t)=\frac{1}{2}\int\limits_{0}^{t}dt_1\int\limits_{0}^{t_1}dt_2\,[Q(t_1), Q(t_2)],$$
$$\Omega_3(t)=\frac{1}{6}\int\limits_{0}^{t}dt_1\int\limits_{0}^{t_1}dt_2\int\limits_{0}^{t_3}dt_3\,\left([Q(t_1),[Q(t_2), Q(t_3)]]+[Q(t_3),[Q(t_2),Q(t_1)]]\right).$$
Square brackets $[A,B]=AB-BA$ are the matrix commutator of $A$ and $B$.

If we represent a matrix $Q(t)$ at the center of the integration interval $t/2$ by Taylor series and keep the main terms with respect to the small parameter $t$, then we get
\begin{equation}\label{mag4}
\Omega_1=tQ+\frac{t^3}{24}Q^{(2)}+O(t^5),\quad \Omega_2=\frac{t^3}{12}\left[Q^{(1)}Q-QQ^{(1)}\right]+O(t^5),\quad
\Omega_3=O(t^5).
\end{equation}
For $t=\tau$ formulas (\ref{mag4}) coincide with (\ref{exp4})-(\ref{power4}). For the Schr\"odinger equation with a time-dependent operator, decomposition was obtained in a series of papers~\cite{Puzynin1999, Puzynin2000}, but the authors decomposed the matrix~$Q$ in the Galerkin series and did not use difference schemes to find derivatives of $Q$.

\subsection{Triple-exponential scheme}

Equation (\ref{T4se1}) can be continued in a different way:
\begin{equation}\label{quasiexp}
T=e^{\tau Q}+\frac{\tau^2}{12}\left[Q^{(1)}e^{\tau Q}-e^{\tau Q}Q^{(1)}\right]
+\frac{\tau^3}{48}\left[e^{\tau Q}Q^{(2)}+Q^{(2)}e^{\tau Q}\right]+O(\tau^5).
\end{equation}
We can continue the scheme (\ref{quasiexp}) to the triple-exponential fourth-order scheme (TES4)
\begin{equation}\label{tripleexp}
T=\exp\left\{\frac{\tau^2}{12}Q^{(1)}+\frac{\tau^3}{48}Q^{(2)}\right\}\exp\left\{\tau Q\right\}\exp\left\{-\frac{\tau^2}{12}Q^{(1)}+\frac{\tau^3}{48}Q^{(2)}\right\},
\end{equation}
which contains a spectral parameter $\zeta $ only at the exponential $e^{\tau {\it Q}} $. The exponential can be split~\cite{Prins2018a}, so the fast techniques (FNFT) can be applied to this scheme.

If $Q$ is skew-Hermitian, then $Q^{(1)}$ and $Q^{(2)}$ are also skew-Hermitian. All exponentials in (\ref{tripleexp}) are unitary matricies and, therefore, they preserve the quadratic invariant. The Maclaurin series of $T$ (\ref{tripleexp}) in $\tau $ gives exactly the decomposition of fourth order schemes for the symmetric case.

\section{Numerical experiments}

\subsection{Numerical algorithms}
We solve the system~(\ref{psit}) on the uniform grid $t_n=-L+\tau n$ with a step~$\tau$ on the interval $[-L,L]$, $L = 30$ unless otherwise stated. If the total number of points is~$2M+1$, then the grid step is $\tau=L/M$.
We replace the original system (\ref{psit}) on the interval $(t_n-\frac{\tau}{2},t_n+\frac{\tau}{2})$ with an approximate system with constant coefficients
\begin{equation}\label{T0}
\Psi_{n+\frac{1}{2}}=T \Psi_{n-\frac{1}{2}}.
\end{equation}
The transition matrix $T$ from the layer $n-\frac{1}{2}$ to the layer $n+\frac{1}{2}$ can be found using different numerical algorithms. Here we compared the numerical results for two new schemes presented above: the exponential scheme ES4 (\ref{exp4}) and triple-exponential scheme TES4 (\ref{tripleexp}). Then we tried the fourth-order conservative transformed scheme (CT4) with the transition operator
\begin{equation}\label{CT4}
T=e^{\frac{\tau}{2}Q_n}\left[I-\frac{\tau}{48}\left(M_{n+1}+M_{n-1}\right)\right]^{-1}
\left[I+\frac{\tau}{48}\left(M_{n+1}+M_{n-1}\right)\right]e^{\frac{\tau}{2}Q_n}
\end{equation}
where
$$M_{n+1}=e^{-\tau Q_n}\left(Q_{n+1}-Q_n\right)e^{\tau Q_n},\quad M_{n-1}=e^{\tau Q_n}\left(Q_{n-1}-Q_n\right)e^{-\tau Q_n}.$$
The CT4 scheme was introduced recently in~\cite{Medvedev2019_OL}. Here we present new and more detailed numerical results for this scheme.

Among the well known algorithms we chose the Boffetta-Osborn second-order scheme (BO)~\cite{Boffetta1992a} and the Runge-Kutta fourth-order algorithm (RK4).
Following~\cite{Burtsev1998} for RK4 scheme we solve the system for the envelope $\chi_1 = \phi_1 e^{i\zeta t}$, $\chi_2 = \phi_2 e^{-i\zeta t}$.
Unlike the above schemes, RK4 does not require computing the transition matrix $T$. Note also that the conventional Runge-Kutta algorithm uses half-steps in its description. Here we set this half-step equal~$\tau$, where $\tau$ is a grid step for the potential~$q(t)$.

The spectral data are finally defined by
\begin{equation}\label{ab_compute}
a(\zeta)= \psi_1(L-\tau/2,\zeta)\,e^{i\zeta (L-\tau/2)},\quad
b(\zeta)= \psi_2(L-\tau/2,\zeta)\,e^{-i\zeta (L-\tau/2)}.
\end{equation}

\subsection{Computation of the derivative of $a(\zeta)$}
To calculate the phase coefficients $r_k$ we need to find the derivatives
\begin{equation}
\frac{da}{d\zeta} = \frac{d\psi_1}{d\zeta}e^{i\zeta (L-\tau/2)} + i(L-\tau/2) a(\zeta).
\end{equation}
From (\ref{T0}) we get
\begin{equation}\label{dpsi}
\frac{d}{d\zeta}\Psi_{n+\frac{1}{2}} = T'_\zeta \Psi_{n-\frac{1}{2}} + T\frac{d}{d\zeta}\Psi_{n-\frac{1}{2}},
\end{equation}
where initial value is defined from (\ref{psi0})
\begin{equation}\label{dpsi0}
\frac{d}{d\zeta}\Psi(-L-\tau/2,\zeta)=\left(\begin{array}{c}-i (-L-\tau/2)\psi_1(-L-\tau/2,\zeta)\\0\end{array}\right).
\end{equation}

For the exponential scheme ES4 (\ref{exp4}) the transition matrix $T$ can be represented as $T=\exp(A)$ and calculated using Pauli matrices (\ref{matrix_exp}) from the appendix. Hence, we can find the derivative
\begin{equation}
T'_{\mbox{(ES4)}}=e^{a_0}
\left(
\begin{array}{cc}
c'+s'a_3 +sa_3' & s'a_1 + sa_1' - is'a_2 - isa_2'\\
s'a_1 + sa_1' + is'a_2 + isa_2' & c'-s'a_3 - sa_3'
\end{array}
\right),
\end{equation}
where
$$
\begin{array}{l}
\displaystyle
c = \cos(\omega),\quad s = \frac{\sin(\omega)}{\omega},\quad \omega=\sqrt{-a_{1}^{2} -a_{2}^{2} -a_{3}^{2} },\\[2mm]
\displaystyle
c' = -\sin(\omega)\omega',\quad s' = \frac{\omega'}{\omega}(c-s),\quad \omega' = -\frac{1}{\omega}(a_1 a_1' + a_2 a_2' + a_3 a_3'),\\[2mm]
\displaystyle
a_1' = i\frac{\tau^2}{24}(d_{12} - d_{21}),\quad a_2' = -\frac{\tau^2}{24}(d_{12} + d_{21}),\quad a_3' = -i\tau,\\[2mm]
\displaystyle
d_{12} = q_{n+1} - q_{n-1},\quad d_{21} = -\sigma(q_{n+1}^* - q_{n-1}^*).
\end{array}
$$

For triple-exponential scheme TES4 (\ref{tripleexp}) we obtain the derivative of the transition matrix~$T$
\begin{equation}
T'_{\mbox{(TES4)}}=\exp\left\{\frac{\tau^2}{12}Q^{(1)}+\frac{\tau^3}{48}Q^{(2)}\right\}
D
\exp\left\{-\frac{\tau^2}{12}Q^{(1)}+\frac{\tau^3}{48}Q^{(2)}\right\},
\end{equation}
where $D = \left(e^{\tau Q_n}\right)'$ can be found from
\begin{equation}\label{dexp}
\left(e^{\tau Q_n}\right)'=-\frac{\tau\zeta}{\omega_n}\sin(\omega_n\tau)I+
\frac{\zeta}{\omega_n^3}
\left[\tau\omega_n\cos(\omega_n\tau)-\sin(\omega_n\tau)\right]Q_n
-i\frac{\sin(\omega_n\tau)}{\omega_n}\sigma_3,
\end{equation}
where $\sigma_3$ is a Pauli matrix (\ref{pauli}).

For CT4 scheme the transition matrix $T$ (\ref{CT4}) can be represented as a fraction $T=A^{-1}B$, therefore
\begin{equation}
T'_{\mbox{(CT4)}}=\left(A^{-1}B\right)'=A^{-1}Q'-A^{-1}A'A^{-1}B,
\end{equation}
where
$$A = \left[I-\frac{\tau}{48}\left(M_{n+1}+M_{n-1}\right)\right]e^{\frac{\tau}{2}Q_n} = A_1 e^{\frac{\tau}{2}Q_n},$$
$$B = \left[I+\frac{\tau}{48}\left(M_{n+1}+M_{n-1}\right)\right]e^{\frac{3\tau}{2}Q_n} = B_1 e^{\frac{3\tau}{2}}.$$
$$A'=-\frac{\tau}{48}\left(\left(e^{2\tau Q_n}\right)'\left(Q_{n-1}-Q_n\right)e^{-2\tau Q_n}+e^{2\tau Q_n}\left(Q_{n-1}-Q_n\right)\left(e^{-2\tau Q_n}\right)'\right)e^{\frac{\tau}{2}Q_n}+A_1\left(e^{\frac{\tau}{2}Q_n}\right)'.$$
$$B'=
\frac{\tau}{48}\left(\left(e^{2\tau Q_n}\right)'\left(Q_{n-1}-Q_n\right)e^{-2\tau Q_n}+e^{2\tau Q_n}\left(Q_{n-1}-Q_n\right)\left(e^{-2\tau Q_n}\right)'\right)e^{\frac{3\tau}{2}Q_n}+B_1\left(e^{\frac{3\tau}{2}Q_n}\right)'.$$
We find the derivative of matrix exponential using (\ref{dexp}).

For BO scheme the derivative of the transition matrix can be found in~\cite{Boffetta1992a}.
For RK4 scheme the derivative of $a(\zeta)$ can be computed using Romberg algorithm~\cite{Burtsev1998, EngelnMullges1996}.

\subsection{Model signals}

We considered a model signal in the form of a chirped hyperbolic secant
\begin{equation}\label{Chirped}
q(t) = A[\mbox{sech}(t)]^{1+iC}.
\end{equation}
For $C = 0$ it is a well-known Satsuma-Yajima signal. The detailed numerical results for this potential are presented in~\cite{Vasylchenkova2019a}.

Here we consider two test potentials: $A = 5.25$, $C = 0$ for anomalous dispersion $\sigma = 1$ and $A = 5.2$, $C = 4$ for both anomalous and normal dispersion $\sigma = \pm 1$.

The analytical expressions of the spectral data of the potential (\ref{Chirped}) for anomalous dispersion are presented in~\cite{Grunbaum1989}. But they can be obtained similarly for normal dispersion ($\sigma = -1$). Here we present general formulas using the Euler Gamma function~$\Gamma$:
\begin{equation}
\begin{array}{l}\label{exact_ab}
\displaystyle
a(\xi)=\frac{\Gamma[1/2-i(\xi + C/2)]\, \Gamma [1/2-i(\xi - C/2)]}
{\Gamma[1/2-i\xi - D]\, \Gamma[1/2-i\xi + D]},\\[5mm]
\displaystyle
b(\xi)=\frac{i}{2^{iC}A} \frac{\Gamma[1/2-i(\xi + C/2)]\, \Gamma [1/2-i(\xi - C/2)]}
{\Gamma[-iC/2 - D]\, \Gamma[-iC/2 + D]},
\quad D = \sqrt{\sigma A^2 - C^2/4}.
\end{array}
\end{equation}

The discrete spectrum $\zeta_k$, $k=\overline{0,K-1}$ is determined by the zeros of the coefficient $a(\zeta)$ and exists only for anomalous dispersion $(\sigma=1)$:
\begin{equation}\label{exact_dzeta}
\zeta_k = i\bigl(\sqrt{A^2 - C^2/4} - 1/2 - k\bigr),\quad k = 0,\ldots,[\sqrt{A^2 - C^2/4} - 1/2],
\end{equation}
where square brackets denote the integer part of the expression.

To compute the phase coefficients~$r_k$ we need to know the derivative~$a'(\zeta)$ only at points~$\zeta_k$ of the discrete spectrum. Let us write the coefficient~$a(\zeta)$ in the following form:
$$a(\zeta)=\frac{f(\zeta)}{\Gamma[1/2-i\zeta - D]},\quad\mbox{where}\quad f(\zeta)=\frac{\Gamma[1/2-i(\zeta + C/2)]\, \Gamma[1/2-i(\zeta - C/2)]}
{\Gamma[1/2-i\zeta + D]}.$$
The function $f(\zeta)$ and its derivative $f'(\zeta)$ have no singularities, so we have at the points of the discrete spectrum:
\begin{equation}
\left.a'(\zeta)\right|_{\zeta=\zeta_k}=-if(\zeta_k)\varphi_k,
\end{equation}
where the function $\varphi_k$ is defined at the points of the discrete spectrum by a recurrence relation
\begin{equation}
\varphi_{k+1}=-(k+1)\varphi_k,\quad \varphi_0=1.
\end{equation}
Thus we have a formula to compute the phase coefficients
\begin{equation}\label{exact_residual}
r_k=\left.\frac{b(\zeta)}{a'(\zeta)}\right|_{\zeta=\zeta_k}=\left.\frac{b(\zeta)}{f(\zeta)}\right|_{\zeta=\zeta_k}\frac{i}{\varphi_k}.
\end{equation}

To calculate energy of the discrete and continuous spectra $E_d$, $E_c$, we use the formula~(\ref{Parseval}): $C_0 = E = E_c + E_d$. Full energy of the potential (\ref{Chirped}) is easily computed as $E = 2A^2$.

Let us denote $K=\left[\sqrt{A^2-C^2/4}+1/2\right]$ as an integer
part of the expression in square brackets and $\delta=\left\{\sqrt{A^2-C^2/4}+1/2\right\}$ as its fractional part. Then the discrete spectrum energy is
\begin{equation}
E_d=4\sum\limits_{k=0}^{K-1}\,\eta_k=2\left(K+\delta-1/2\right)^2-2\left(\delta-1/2\right)^2.
\end{equation}
and the continuous spectrum energy is
\begin{equation}
E_c=E-E_d=2\left(C^2/4+\left(\delta-1/2\right)^2\right).
\end{equation}

\subsection{Approximation order}
The following formula was used to calculate the approximation order $m$:
\begin{equation}\label{order}
m=\log_{\frac{\tau_1}{\tau_2}}\frac{\left\Vert\tilde{\Psi}_1(L)\right\Vert_2}
{\left\Vert\tilde{\Psi}_2(L)\right\Vert_2}=
\frac{\log_2\frac{\left\Vert\tilde{\Psi}_1(L)\right\Vert_2}
	{\left\Vert\tilde{\Psi}_2(L)\right\Vert_2}}{\log_2(\tau_1/\tau_2)},
\end{equation}
where $\tau_i$, $i=1$, $2$ are the steps of computational grids for two calculations with one spectral parameter $\zeta$ and $\tau_1 > \tau_2$, $\tilde{\Psi}_i(L)$ is a deviation of the calculated value $\Psi_i(L)$ from the exact analytical value $\bar{\Psi}_i(L)$ at the boundary point $t=L$. The calculations were carried out for different $p$-norms and showed close values for the approximation orders. However, for the Euclidean $2$-norm, the graphics were the smoothest.

\begin{figure}
  \centering
  \fbox{\includegraphics{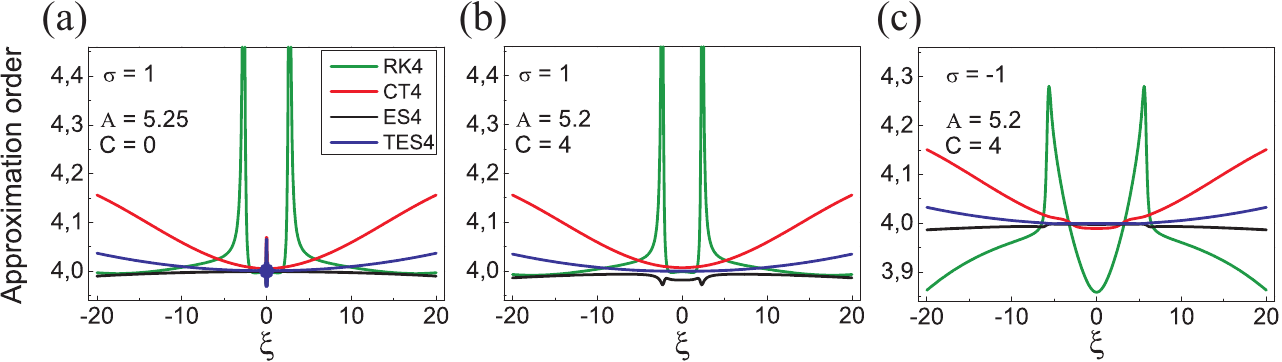}}
  \caption{The approximation order with respect to the spectral parameter $\xi$.}
 	\label{fig:1}
\end{figure}

Figure~\ref{fig:1} confirms the approximation order~$m = 4$ of the schemes with respect to a spectral parameter $\xi\in [-20, 20]$. Each line was calculated by the formula~(\ref{order}) using two embedded grids with a doubled grid step $\tau = L/M$, $L = 30$, where coarse and fine grids were defined by $M = 2^{10}$ and $M=2^{11}$. Let us remind that the total number of points in the whole domain~$[-L, L]$ is $2M +1$.

\subsection{Formulas for errors}
We present the numerical errors of calculating the spectral data for continuous and discrete spectrum. To find the calculation errors of the continuous spectrum energy~$E_c$, residuals~$r_k$, and the coefficients $a(\zeta)$, $b(\zeta)$ at fixed~$\zeta$ we use formula
\begin{equation}\label{error}
\mbox{error}[\phi]=\frac{|\phi^{comp} - \phi^{exact}|}{|\phi_0|},\quad
\phi_0 =
\begin{cases}
\phi^{exact}, \mbox{ if } |\phi^{exact}| > 1\\
1, \mbox{ otherwise},
\end{cases}
\end{equation}
where $\phi$ can represent $E_c$, $r_k$, $a(\zeta)$ or $b(\zeta)$ at fixed $\zeta$.

For the continuous spectrum we calculate the mean squared error
\begin{equation}\label{MSE}
MSE[\phi]=\frac{1}{N}\sum_{j=1}^{N}\frac{|\phi^{comp}(\xi_j) - \phi^{exact}(\xi_j)|^2}{|\phi_0(\xi_j)|^2},\quad
\phi_0 =
\begin{cases}
\phi^{exact}(\xi_j), \mbox{ if } |\phi^{exact}(\xi_j)| > 1\\
1, \mbox{ otherwise},
\end{cases}
\end{equation}
where $\phi$ can represent $a(\xi)$ or $b(\xi)$. Here we suppose the spectral parameter $\xi\in [-20, 20]$ with the total number of points $N = 1025$.

\subsection{Numerical results for continuous spectrum}

\begin{figure}
  \centering
  \fbox{\includegraphics{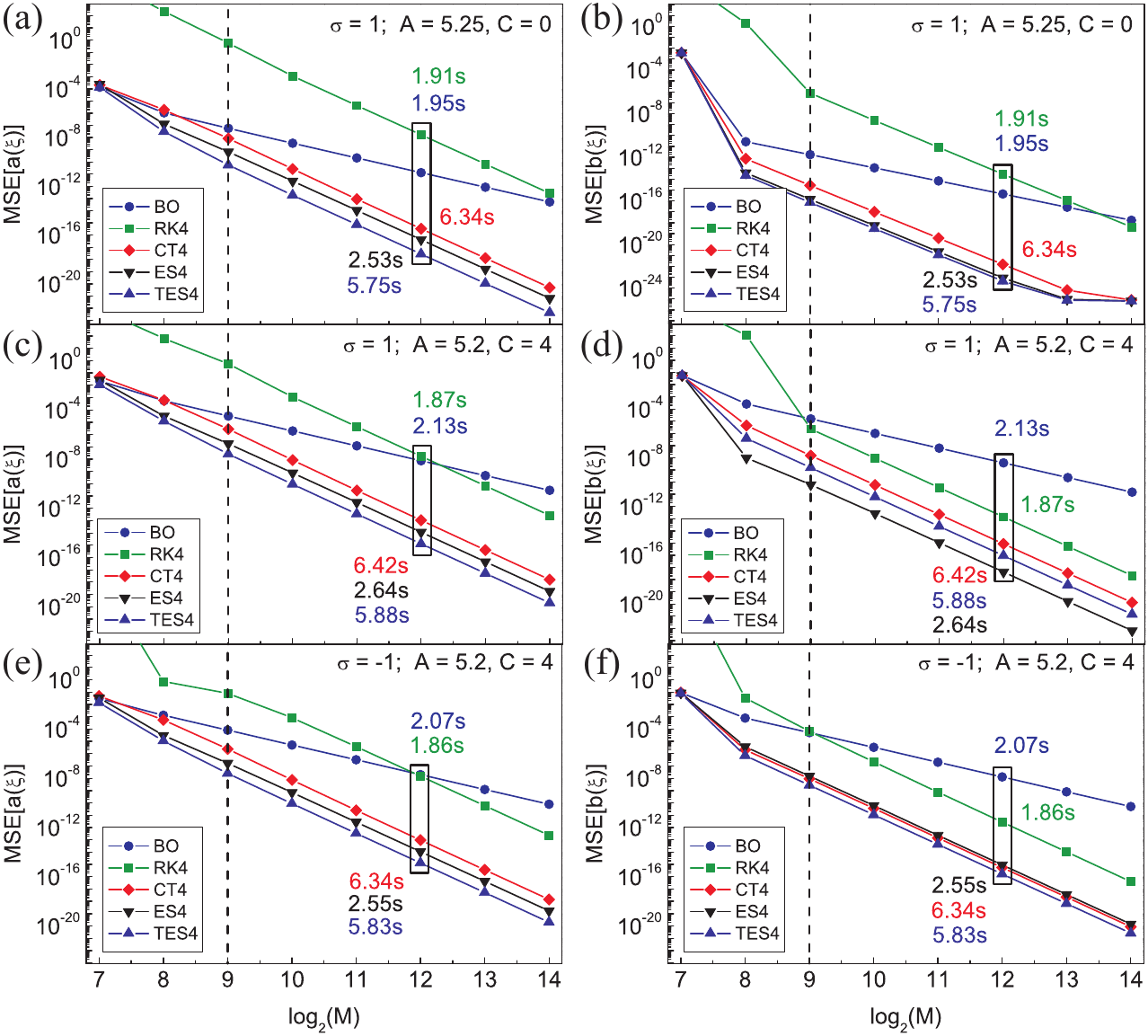}}
  \caption{Continuous spectrum mean squared errors (\ref{MSE}) of $a(\xi)$ and $b(\xi)$.}
 	\label{fig:2}
\end{figure}

\begin{figure}
  \centering
  \fbox{\includegraphics{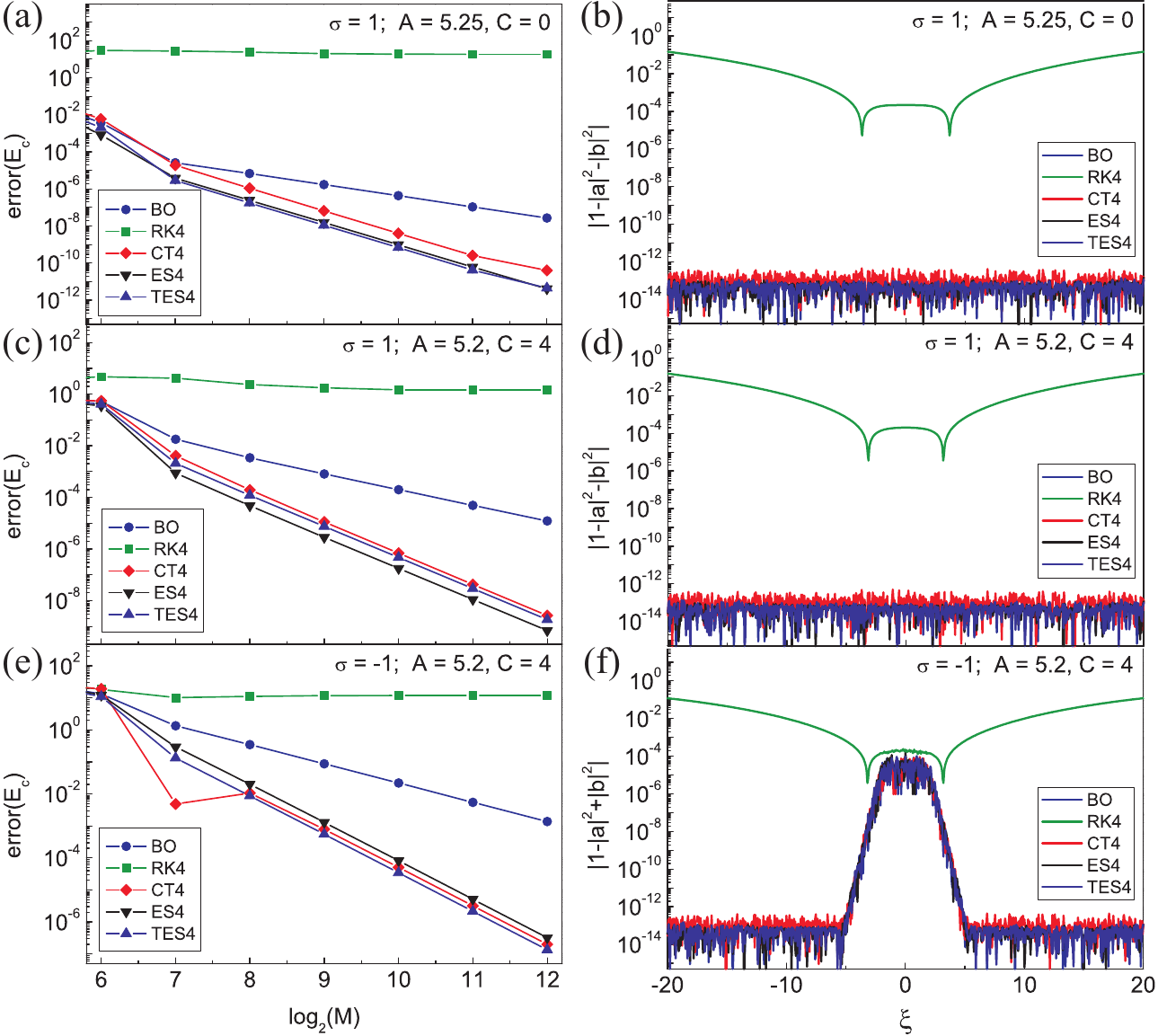}}
  \caption{(a, c, e) Errors (\ref{error}) of continuous spectrum energy.
		(b, d, f) Absolute errors of quadratic invariants.}
 	\label{fig:3}
\end{figure}

Figures~\ref{fig:2},~\ref{fig:3} present the continuous spectrum errors for the potential~(\ref{Chirped}) with
two sets of parameters: $A = 5.25$, $C = 0$ for anomalous dispersion $\sigma = 1$ and $A = 5.2$, $C = 4$ for both anomalous and normal dispersion $\sigma = \pm 1$.

Figure~\ref{fig:2} shows the mean squared error~(\ref{MSE}) of the coefficients $a(\xi)$ and $b(\xi)$ with respect to the number of grid nodes $M$.
The total number of points in the whole domain $[-L,L]$ is~$2M+1$. Dashed vertical lines mark the minimum number of grid nodes~$M_{\min}$ that guarantee a good approximation~\cite{Medvedev2019_OL}. Actually, when calculating the continuous spectrum, it is necessary to choose a time step~$\tau=L/M$ to describe correctly the fastest oscillations. For a fixed value of~$\xi$, the local frequency $\omega(t;\xi)=\sqrt{\xi^2+|q(t)|^2}$ of the system~(\ref{psit}) varies from $\omega_{\min}=|\xi|$ to $\omega_{\max}=\sqrt{\xi^2+q_{\max}^2}$, where $q_{\max}=\max\limits_t|q(t)|$ is the maximum absolute value of the potential $q(t)$. Therefore, step $\tau$ cannot be arbitrary. In order to describe the most rapid oscillations, it is necessary to have at least 4-time steps for the oscillation period, so the inequality must be satisfied:
$$4\tau=4\frac{L}{M}\leq \frac{2\pi}{\omega_{\max}}.$$
Therefore, any difference schemes will approximate the solutions of the original continuous system~(\ref{psit}) if the inequality is fulfilled for the number of points $M\geq M_{\min}=2\,L\,\omega_{\max}/\pi$.

Figure~\ref{fig:2} also demonstrates a comparison of the computational time. One can see that ES4 scheme shows the best accuracy with a maximum speed.

Figure~\ref{fig:3} shows how the numerical schemes conserve energy.
The numerical errors (\ref{error}) for the continuous spectrum energy (\ref{Econt}) are compared in Fig.~\ref{fig:3} (a, c, e).
To calculate the continuous spectrum energy it is important to define the size of the spectral domain~$L_{\xi}$ and the corresponding grid step~$d\xi$~\cite{Medvedev2019_OL}. According to the conventional discrete Fourier transform, we take the same number of points $N_{\xi} = N$ in the spectral domain and define a spectral step as $d\xi = \pi/(2L)$. So the size of the spectral interval is
\begin{equation}\label{Lxi}
L_{\xi} = \pi/(2\tau).
\end{equation}

Figure~\ref{fig:3} (b, d, f) demonstrates the deviation of the quadratic invariant $H=|a|^2+\sigma |b|^2$ from unit with respect to the real spectral parameter $\xi$.

If the matrix $Q$ is skew-Hermitian, then the matrix $\exp(\tau Q)$ is unitary and the quadratic invariant conserves. For the direct ZSP this corresponds to anomalous dispersion $(\sigma = 1)$ with a real spectral parameter $\zeta=\xi$.

Let us consider a more general system with a matrix $Q=KD$,
where $K(t)$ is anti-Hermitian matrix depending on $t$, $D$ is a constant Hermitian matrix. The system (\ref{psit}) conserves the quadratic value $H=(\Psi^*,D\Psi)$. Indeed, we have this result from the chain of equalities
$$
\frac{d}{dt}\left(\Psi^*,D\Psi\right)=\left(\frac{d\Psi^*}{dt},D\Psi\right)+
\left(\Psi^*,D\frac{d\Psi}{dt}\right)=\left(K^*D^*\Psi^*,D\Psi\right)+\left(\Psi^*,DKD\Psi\right)=
$$
$$=\left({\Psi}^*,(D^*)^T(K^*)^TD{\Psi}\right)+\left({\Psi}^*,DKD{\Psi}\right)=
\left({\Psi}^*,D(K^\dagger+K)D{\Psi}\right)=0.$$

For one-step exponential methods $\Psi_{n+\frac{1}{2}}=e^{R_nD}\Psi_{n-\frac{1}{2}}$, where $R_n$ is skew-Hermitian matrix, the quadratic invariant also conserves. It follows from the chain of equalities
$$
(\Psi_{n+\frac{1}{2}}^*,D\Psi_{n+\frac{1}{2}})=(e^{\tau R_n^*D^*}\Psi_{n-\frac{1}{2}}^*,De^{\tau R_nD}\Psi_{n-\frac{1}{2}})=
(e^{\tau R_n^*D^*}\Psi_{n-\frac{1}{2}}^*,e^{\tau DR_n}D\Psi_{n-\frac{1}{2}})=
$$
$$=(\Psi_{n-\frac{1}{2}}^*,e^{-\tau R_nD}e^{\tau DR_n}D\Psi_{n-\frac{1}{2}})=(\Psi_{n-\frac{1}{2}}^*,D\Psi_{n-\frac{1}{2}}).$$
Here we used the formula $De^{\tau R_nD}=e^{\tau DR_n}D$, because for any natural $p$ the equality is valid: $D(R_nD)^p=(DR_n)^pD$.

From this result follows, that Boffetta-Osborn scheme (\ref{T0}) and the exponential scheme (\ref{exp4}) are conservative for normal and anomalous dispersion $\sigma=\pm 1$. Similarly the scheme (\ref{CT4}) also conserves the quadratic invariant, because it is the function of $R_nD$.
\begin{figure}
  \centering
  \fbox{\includegraphics{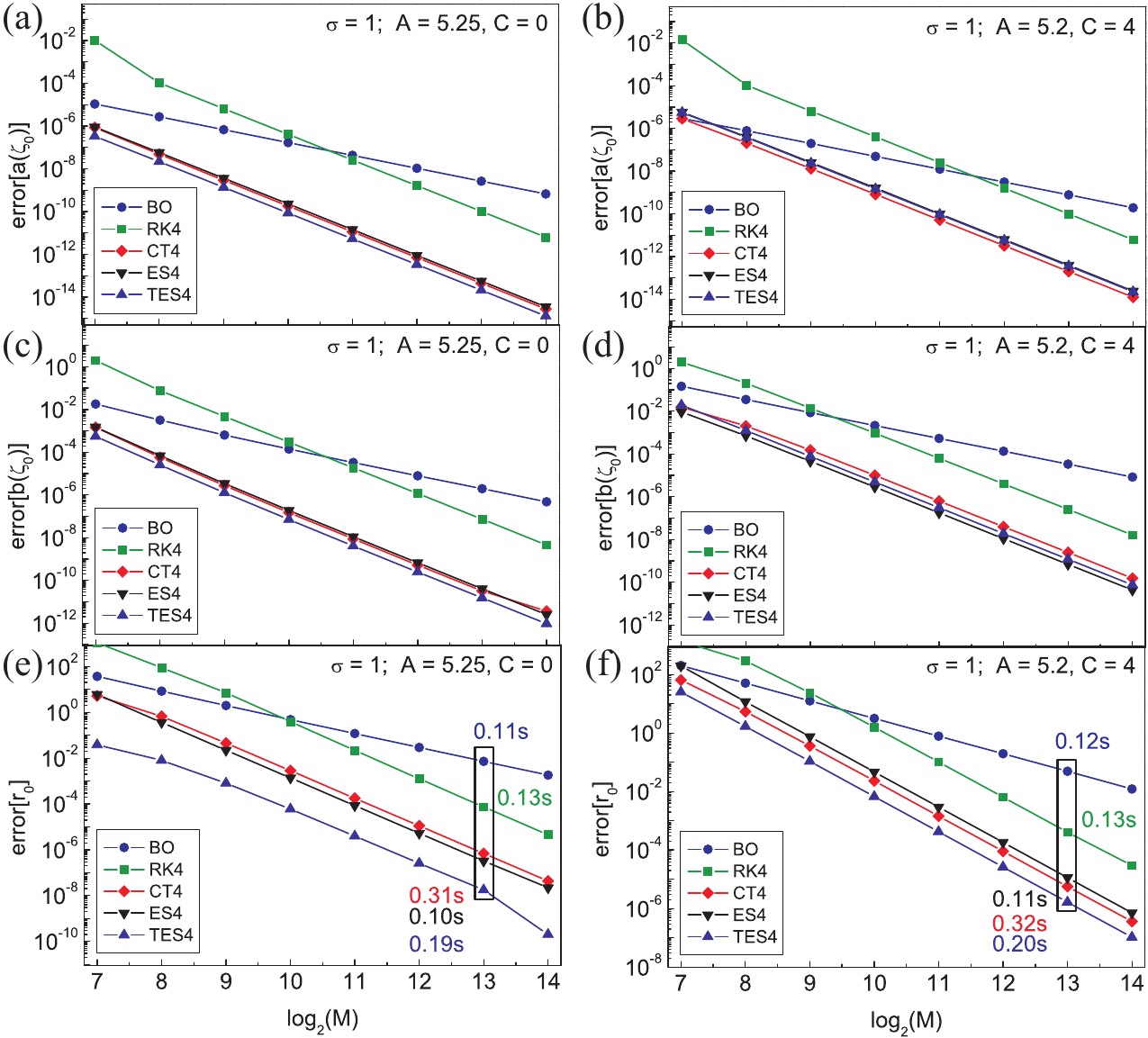}}
  \caption{Discrete spectrum errors for the maximum eigenvalue $\zeta_0$.}
 	\label{fig:4}
\end{figure}

Figure~\ref{fig:3} confirms that RK4 scheme does not conserve the continuous spectrum energy and quadratic invariant for the real spectral parameters.

Figure~\ref{fig:3} (f) corresponds to the case of normal dispersion, therefore in the center of the spectral interval the parameters $a(\xi)$ and $b(\xi)$ have large values. This leads to the higher computational error in this zone. At the same time the quadratic invariant in this case equally conserves for all schemes considered here. However, RK4 scheme again shows the worst results at the edges of the spectral interval.

\subsection{Numerical results for discrete spectrum}

\begin{figure}
  \centering
  \fbox{\includegraphics{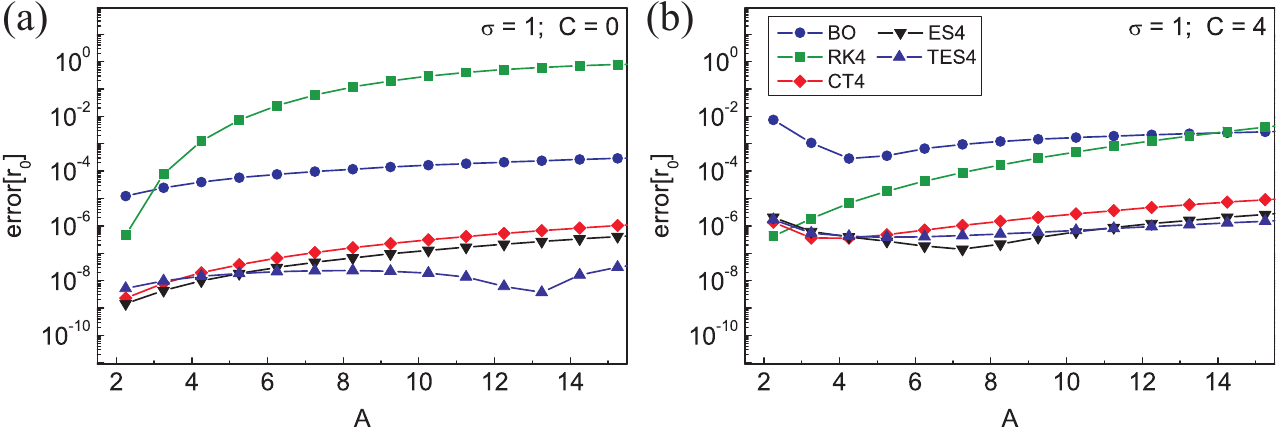}}
  \caption{The errors (\ref{error}) of phase coefficients for the maximum eigenvalue $\zeta_0$ (\ref{exact_dzeta}).}
 	\label{fig:5}
\end{figure}

Figures~\ref{fig:4},~\ref{fig:5} present the discrete spectrum errors (\ref{error}). The parameters $a(\zeta_k)$, $b(\zeta_k)$ and $r_k$ were computed for the analytically known eigenvalues (\ref{exact_dzeta}). Here we did not use any numerical algorithm to find the eigenvalue but computed spectral data at the exact point $\zeta = \zeta_k$ right away. It was made intentionally to estimate the error of the scheme itself and to avoid the influence of the other numerical algorithm errors.

There are well known problems with the computation of the coefficient $b(\zeta_k)$. We used the bi-directional algorithm~\cite{7486016} to find $b(\zeta_k)$ by the formula (\ref{bidir}).

Figure~\ref{fig:5} presents the errors (\ref{error}) of computing the phase coefficients for the maximum eigenvalue $\zeta_0$ (\ref{exact_dzeta}) with respect to the amplitude $A$ of the potential (\ref{Chirped}). Here $L = 20$, $M = 2^{11}$.

\section{Conclusion}
Two new forth-order exponential schemes for the numerical solution of the direct Zakharov-Shabat problem were presented and compared with the known ones. The ES4 scheme demonstrates the excellent computational speed and accuracy, but has difficulties with the direct application of the fast algorithms.

The TES4 scheme shows the same accuracy, but it requires about $2$ times longer to calculate. However the main advantage of TES4 scheme is that fast algorithms can be applied to it.

The CT4 scheme also shows good accuracy comparable with two exponential schemes mentioned above, but it works about $2.5-3$ times longer than ES4. The CT4 scheme does not allow the direct application of the fast algorithms. But it is possible after exponential approximation applied to the scheme.

All these schemes have an advantage over RK4 because they conserve the energy for the continuous spectrum parameter.

\appendix
\section*{Appendix}
The Pauli matrices can be used to calculate the matrix exponential in (\ref{exp4}). A $2\times 2$ complex matrix $A$ can be presented as
$A=a_{0}\sigma_{0} +a_{1}\sigma_{1} +a_{2}\sigma_{2} +a_{3} \sigma_{3}$,
where $\sigma_{j} $ are Pauli matrices:
\begin{equation}\label{pauli}
\sigma_0=\left(\begin{array}{cc}1&0\\0&1\end{array}\right),\quad \sigma_1=\left(\begin{array}{cc}0&1\\1&0\end{array}\right),\quad
\sigma_2=\left(\begin{array}{cc}0&-i\\i&0\end{array}\right),\quad
\sigma_3=\left(\begin{array}{cc}1&0\\0&-1\end{array}\right)
\end{equation}
Then the matrix exponential can be found as
\begin{equation}\label{matrix_exp}
e^{A}=e^{a_0}\left[c\sigma_0+
is\left(a_1\sigma_1+a_2\sigma_2+a_3\sigma_3\right)\right]=
e^{a_0}
\left(
\begin{array}{cc}
c+sa_3 & s(a_1 - i a_2)\\
s(a_1 + i a_2) & c-sa_3
\end{array}
\right),
\end{equation}
where $\omega=\sqrt{-a_{1}^{2} -a_{2}^{2} -a_{3}^{2} } $, $c = \cos(\omega)$, $\displaystyle s = \frac{\sin(\omega)}{\omega}$.

\section*{Funding}
Russian Science Foundation (RSF) (17-72-30006).

\bibliographystyle{unsrt}

\begin{thebibliography}{10}

\bibitem{ZakharovShabat1972}
V.~E. Zakharov and A.~B. Shabat.
\newblock {Exact Theory of Two-Dimensional Self-Focusing and One-Dimensional
  Self-Modulation of Waves in Non-Linear Media}.
\newblock {\em Journal of Experimental and Theoretical Physics}, 34(1):62--69,
  1972.

\bibitem{Yousefi2014II}
Mansoor~I Yousefi and Frank~R Kschischang.
\newblock {Information Transmission Using the Nonlinear Fourier Transform, Part
  II: Numerical Methods}.
\newblock {\em IEEE Transactions on Information Theory}, 60(7):4329--4345,
  2014.

\bibitem{Turitsyn2017Optica}
Sergei~K. Turitsyn, Jaroslaw~E. Prilepsky, Son~Thai Le, Sander Wahls, Leonid~L.
  Frumin, Morteza Kamalian, and Stanislav~A. Derevyanko.
\newblock {Nonlinear Fourier transform for optical data processing and
  transmission: advances and perspectives}.
\newblock {\em Optica}, 4(3):307, 2017.

\bibitem{Vasylchenkova2019a}
A~Vasylchenkova, J.E. Prilepsky, D~Shepelsky, and A~Chattopadhyay.
\newblock {Direct nonlinear Fourier transform algorithms for the computation of
  solitonic spectra in focusing nonlinear Schr{\"{o}}dinger equation}.
\newblock {\em Communications in Nonlinear Science and Numerical Simulation},
  68:347--371, 3 2019.

\bibitem{Boffetta1992a}
G.~Boffetta and A.R Osborne.
\newblock {Computation of the direct scattering transform for the nonlinear
  Schroedinger equation}.
\newblock {\em Journal of Computational Physics}, 102(2):252--264, 10 1992.

\bibitem{Medvedev2019_OL}
Sergey Medvedev, Irina Vaseva, Igor Chekhovskoy, and Mikhail Fedoruk.
\newblock {Numerical algorithm with fourth-order accuracy for the direct
  Zakharov-Shabat problem}.
\newblock {\em Optics Letters}, 44(9):2264, 2019.

\bibitem{Ablowitz1981}
Mark~J. Ablowitz and Harvey Segur.
\newblock {\em {Solitons and the Inverse Scattering Transform}}.
\newblock Society for Industrial and Applied Mathematics, Philadelphia, 1981.

\bibitem{Dahlquist1963}
Germund~G. Dahlquist.
\newblock {A special stability problem for linear multistep methods}.
\newblock {\em BIT}, 3(1):27--43, 3 1963.

\bibitem{Hairer1987}
Ernst Hairer, Syvert~P. N{\o}rsett, and Gerhard Wanner.
\newblock {\em {Solving ordinary differential equations I. nonstiff problems}}.
\newblock Springer-Verlag Berlin Heidelberg, 1987.

\bibitem{Magnus1954}
Wilhelm Magnus.
\newblock {On the exponential solution of differential equations for a linear
  operator}.
\newblock {\em Communications on Pure and Applied Mathematics}, 7(4):649--673,
  1954.

\bibitem{Puzynin1999}
I.V. Puzynin, A.V. Selin, and S.I. Vinitsky.
\newblock {A high-order accuracy method for numerical solving of the
  time-dependent Schr{\"{o}}dinger equation}.
\newblock {\em Computer Physics Communications}, 123(1-3):1--6, 1999.

\bibitem{Puzynin2000}
I.V. Puzynin, A.V. Selin, and S.I. Vinitsky.
\newblock {Magnus-factorized method for numerical solving the time-dependent
  Schr{\"{o}}dinger equation}.
\newblock {\em Computer Physics Communications}, 2000.

\bibitem{Prins2018a}
Peter~J Prins and Sander Wahls.
\newblock {Higher Order Exponential Splittings for the Fast Non-Linear Fourier
  Transform of the Korteweg-De Vries Equation}.
\newblock In {\em ICASSP, IEEE International Conference on Acoustics, Speech
  and Signal Processing - Proceedings}, number~4, pages 4524--4528. IEEE, 2018.

\bibitem{Burtsev1998}
S~Burtsev, R~Camassa, and I~Timofeyev.
\newblock {Numerical Algorithms for the Direct Spectral Transform with
  Applications to Nonlinear Schr{\"{o}}dinger Type Systems}.
\newblock {\em Journal of Computational Physics}, 147(1):166--186, 11 1998.

\bibitem{EngelnMullges1996}
Giesela Engeln-Mullges and Frank Uhlig.
\newblock {\em {Numerical Algorithms with C}}.
\newblock Springer-Verlag Berlin Heidelberg, 1996.

\bibitem{Grunbaum1989}
F.~A. Grunbaum.
\newblock {The scattering problem for a phase-modulated hyperbolic secant
  pulse}.
\newblock {\em Inverse Problems}, 5(3):287--292, 1989.

\bibitem{7486016}
Siddarth Hari and Frank~R. Kschischang.
\newblock {Bi-Directional Algorithm for Computing Discrete Spectral Amplitudes
  in the NFT}.
\newblock {\em Journal of Lightwave Technology}, 34(15):3529--3537, 8 2016.

\end{thebibliography}

\end{document}